\documentclass[12pt]{article} %***
\usepackage[sectionbib]{natbib}
\usepackage{array,epsfig,fancyheadings}
\usepackage[]{hyperref}  %<----modified by Ivan
\usepackage{sectsty, secdot}

\usepackage{amssymb}
\usepackage{subfig}
\usepackage{float}
\usepackage{bm}
\usepackage{enumerate}
\usepackage{algorithm, algorithmic,hyperref}
\usepackage{booktabs,multirow}
\usepackage{diagbox}
\usepackage{mathrsfs}
\usepackage{graphicx,setspace,lscape,longtable}
\usepackage{natbib,epsfig,graphicx}
\usepackage{amsmath,amsthm,amssymb,color}
\usepackage{amsfonts}
\usepackage{siunitx}
\bibpunct{(}{)}{;}{a}{,}{,}
\usepackage{longtable}

\usepackage{siunitx}
\sectionfont{\fontsize{12}{14pt plus.8pt minus .6pt}\selectfont}
\renewcommand{\theequation}{\thesection\arabic{equation}}
\subsectionfont{\fontsize{12}{14pt plus.8pt minus .6pt}\selectfont}
\usepackage{amsmath}
\usepackage{amssymb}
\usepackage{amsfonts}
\usepackage{multirow}
\usepackage{amsthm}

\newtheorem{theorem}{Theorem}
\newtheorem{lemma}{Lemma}
\newtheorem{corollary}{Corollary}
\newtheorem{proposition}{Proposition}
\theoremstyle{definition}
\newtheorem{definition}{Definition}

\def\sxx{\widehat\Sigma_{xx}}
\def\sxy{\widehat\Sigma_{xy}}

\def\lmax{\lambda_{\max}}
\def\lmin{\lambda_{\min}}

\def\ols{{\operatorname{ols}}}

\def\op{\operatorname{op}}

\def\beq{\begin{equation}}
\def\eeq{\end{equation}}
\def\beqs{\begin{equation*}}
\def\eeqs{\end{equation*}}
\def\beqr{\begin{eqnarray}}
\def\eeqr{\end{eqnarray}}
\def\beqrs{\begin{eqnarray*}}
\def\eeqrs{\end{eqnarray*}}
\def\bet{\begin{theorem}}
\def\eet{\end{theorem}}

\def\bec{\begin{corollary}}
\def\eec{\end{corollary}}

\def\bed{\begin{definition}}
\def\eed{\end{definition}}
\def\bel{\begin{lemma}}
\def\eel{\end{lemma}}
\def\bep{\begin{proposition}}
\def\eep{\end{proposition}}
\def\bea{\begin{assumption}}
\def\eea{\end{assumption}}
\def\bg{\begin{figure}[tbph]\begin{center}}
\def\eg{\end{center}\end{figure}}

\def\bc{\begin{center}}
\def\ec{\end{center}}

\def\wt{\widetilde}
\def\wh{\widehat}

\def\lmax{\lambda_{\max}}

\def\mL{\mathcal L}

\def\mR{\mathbb{R}}

\def\mS{\mathcal S}

\def\var{\mbox{var}}

\def\cov{\mbox{cov}}

\def\rank{\operatorname{rank}}
\def\diag{\operatorname{diag}}

\DeclareMathOperator{\tr}{tr}

\def\tr{\operatorname{tr}}

\def\rank{\operatorname{rank}}

\def\bols{\widehat{\beta}_{\operatorname{ols}}}

\def\rank{\operatorname{rank}}
\def\trainid{\mathcal I_{\text{train}}}
\def\testid{\mathcal I_{\text{test}}}
\def\avgerr{\text{AvgErr}}

\def\sub{\operatorname{sub}}

\usepackage[margin=1in]{geometry}

\begin{document}

\pagestyle{plain}	%%%%%%%%%%%%%%%%%%%%%%%%%%%%%%%%%%%%%%%%%%%%%%%%%%%%%%%%%%%%%%%%%%%%%%%%%%%%%%%%%%%%%%%%%%%%%%%%%%%%%%%%%%%%%%%%%%%%%%%%%%%%
	%%%%%%%%%%%%%%%%%%%%%%%%%%%%%%%%%%%%%%%%%%%%%%%%%%%%%%%%%%%%%%%%%%%%%%%%%%%%%%%%%%%%%%%%%%%%%%%%%%%%%%%%%%%%%%%%%%%%%%%%%%%%
	
	\renewcommand{\baselinestretch}{2}
	
	\markright{ \hbox{\footnotesize\rm Statistica Sinica
			%{\footnotesize\bf 24} (201?), 000-000
		}\hfill\\[-13pt]
		\hbox{\footnotesize\rm
			%\href{http://dx.doi.org/10.5705/ss.20??.???}{doi:http://dx.doi.org/10.5705/ss.20??.???}
		}\hfill }
	
	\markboth{\hfill{\footnotesize\rm HANCHAO YAN, FEIFEI WANG ET AL.} \hfill}
	{\hfill {\footnotesize\rm AUXILIARY LEARNING AND ITS STATISTICAL UNDERSTANDING} \hfill}
	
	\renewcommand{\thefootnote}{}
	$\ $\par
	
	%%%%%%%%%%%%%%%%%%%%%%%%%%%%%%%%%%%%%%%%%%%%%%%%%%%%%%%%%%%%%%%%%%%%%%%%%%%%%%%%%%%%%%%%%%%%%%%%%%%%%%%%%%%%%%%%%%%%%%%%%%%%
	
	\fontsize{12}{14pt plus.8pt minus .6pt}\selectfont \vspace{0.8pc}
	\centerline{\large\bf Auxiliary Learning and its Statistical Understanding}
	% \vspace{2pt}
	% \centerline{\large\bf HERE IF A SECOND LINE IS NEEDED}
	\vspace{.3cm}
	\centerline{Hanchao Yan$^{1}$, Feifei Wang$^{2}$\footnote{The corresponding author. Email: feifei.wang@ruc.edu.cn}, Chuanxin Xia$^3$, and Hansheng Wang$^1$    }
	% \vspace{.3cm}
    % \centerline{\it $^1$ Peking University, $^2$ Renmin University of China, $^3$ Ministry of Commerce People's Republic of China}
    \begin{center}
    \it
        $^1$ Peking University; $^2$ Renmin University of China;\\
        $^3$ Ministry of Commerce People's Republic of China.
    \end{center}

 \fontsize{9}{11.5pt plus.8pt minus.6pt}\selectfont
	
	%%%%%%%%%%%%%%%%%%%%%%%%%%%%%%%%%%%%%%%%%%%%%%%%%%%%%%%%%%%%%%%%%%%%%%%%%%%%%%%%%%%%%%%%%%%%%%%%%%%%%%%%%%%%%%%%%%%%%%%%%%%%
	
	\begin{quotation}
		\noindent {\it Abstract:}
        Modern statistical analysis often encounters high-dimensional problems but with a limited sample size. It poses great challenges to traditional statistical estimation methods.
        In this work, we adopt auxiliary learning to solve the estimation problem in high-dimensional settings. We start with the linear regression setup. To improve the statistical efficiency of the parameter estimator for the primary task, we consider several auxiliary tasks, which share the same covariates with the primary task. Then a weighted estimator for the primary task is developed, which is a linear combination of the ordinary least squares estimators of both the primary task and auxiliary tasks. The optimal weight is analytically derived and the statistical properties of the corresponding weighted estimator are studied. We then extend the weighted estimator to generalized linear regression models. Extensive numerical experiments are conducted to verify our theoretical results. Last, a deep learning-related real-data example of smart vending machines is presented for illustration purposes.

		\vspace{9pt}
		\par
	\end{quotation}\par

	\def\thefigure{\arabic{figure}}
	\def\thetable{\arabic{table}}
	
	\renewcommand{\theequation}{\thesection.\arabic{equation}}

	\fontsize{12}{14pt plus.8pt minus .6pt}\selectfont

\newpage	
	\section{INTRODUCTION}
	
	Modern statistical analysis frequently encounters datasets characterized by an extensive number of features, which even surpass the sample size. This is attributed to advancements in state-of-the-art technology, which facilitates more convenient collection of data and computing \citep{gao2022review, jiang2023feature,li2024selective}. For instance, in genetic studies, there are hundreds of thousands of single-nucleotide polymorphisms (SNPs) that serve as potential predictors. However, the sample size is limited due to the expensive cost of sequencing. In deep learning-related studies, millions or even billions of features are extracted from texts, images, and other forms of unstructured data. However, the sample size is also relatively limited due to the high cost of labeling. The fact that the feature dimension is larger than the sample size imposes significant challenges to classical statistical analysis. Consider the linear regression model as an example. To estimate the regression coefficient $\beta$, the popularly used ordinary least squares (OLS) estimator can achieve excellent estimation performance under the setting of fixed $p$ and $N \to \infty$. However, in the high-dimensional setting with $p\to\infty$ as $N\to\infty$, the finite sample performance of the OLS estimator is no longer satisfactory \citep{wainwright2019high}. Then how to improve the estimation accuracy becomes a problem of great interest.

To solve this problem, a variety of statistical methods have been developed, which can be roughly classified into two streams. The first stream is penalized estimation, including LASSO \citep{tibshirani1996regression}, the SCAD \citep{fan2001variable}, group LASSO \citep{Yuan2005}, adaptive LASSO \citep{zou2006adaptive}, prior LASSO \citep{jiang2016variable}, and many others.  The second stream is feature screening. Along this direction, one notable work is the Sure Independence Screening (SIS) method proposed by \cite{fan2008sure}. \citet{wang2009forward} further improved SIS by applying a forward regression algorithm to screen important features and demonstrated its screening consistency property under an ultrahigh-dimensional setup.
The ultrahigh dimensional screening problem but with a low-dimensional factor structure was studied by \citet{wang2012factor}.
\citet{Li2012Feature} proposed a model-free feature screening method based on the distance correlation. Recently, \cite{Wang2021Sequential} developed a sequential term selection method for textual datasets, and \cite{zhu2022feature} discussed this feature screening method in the context of subsampling. More discussions about the statistical methods in these two streams can be found in \citet{fan2020statistical}.

In this work, we try to address this problem from a different perspective. That is, \textbf{auxiliary learning}, which is a popular concept in the field of machine learning \citep{kung-etal-2021-efficient, 9643357, Chen2022Module}.
In auxiliary learning, the research of main interest is often referred to as the primary task. In addition to the primary task, several auxiliary tasks are also considered to assist the primary task in learning better representations or improving the overall performance. 
To optimize the performance for the primary task, previous researches about auxiliary learning often adopt weighted loss functions to leverage the importance of auxiliary tasks. For example, \citet{liebel2018auxiliary} directly treated the weights of different tasks as parameters to be learned with certain regularization.
\citet{navon2020auxiliary} adopted a separate network structure to combine
all losses into one single objective function, so that nonlinear interactions between different tasks can be learned.
\citet{chen2022auxiliary} also adopted the approach of using separate network. However, they focused on evaluating the relevance of data features to different tasks, which were then combined to form the weighted loss function. 

Note that, auxiliary learning is a special practice of multi-task learning. However, different from auxiliary learning, in multi-task learning, the model is trained on multiple tasks simultaneously to improve the performance of each task by leveraging shared representations \citep{vandenhende2021multi}. To achieve representation sharing, a commonly used technique is to use hard-parameter sharing \citep{kokkinos2017ubernet, chen2018gradnorm,kendall2018multi} or soft-parameter sharing \citep{misra2016cross, ruder2019latent}. Both techniques mainly focus on the shared layers between the tasks to improve the performance of the specific tasks.

Auxiliary learning has been widely adopted in practice \citep{kanezaki2018rotationnet, howard2018universal}.
Take the real application of smart vending machines in Section 4 as an example. The primary task here is to accurately identify each product of bottled water via the shopping images. It is a multi-class classification problem. Preliminary findings revealed that solely classifying this task yielded suboptimal results. To enhance the classification accuracy, several auxiliary tasks can be considered. For example, we observed that factors such as the cap color, the liquid color, the label shape, whether there are horizontal ridges or vertical ridges contribute to identifying the product of bottled water. Thus identifying these factors can be regarded as auxiliary tasks to enhance the learning process of the primary tasks.

In this work, we aim to study auxiliary learning from a statistical perspective. We start with the linear regression setup for its analytical simplification. However, both our numerical and real data experiments suggest that the theoretical insights obtained here can be readily applied to more sophisticated models (e.g., logistic regression and deep learning). Let $(X_i,Y_i^{(0)})$ be the observation collected from the $i$-th subject with $1\leq i \leq N$. Here $Y_i^{(0)}$ is a continuous response of the primary task and $X_i\in \mR^p$ is the associated covariate. To model the regression relationship between $X_i$ and $Y_i^{(0)}$, a standard linear regression model can be assumed to be $Y_i^{(0)}=X_i^\top\beta^{(0)}+\varepsilon_i^{(0)}$, where $\beta^{(0)}$ is the primary regression coefficient and $\varepsilon_i^{(0)}$ is the independent error term. Except for $Y_i^{(0)}$, we also consider $K$ auxiliary responses $Y_i^{(k)}$ with $1\leq k \leq K$. Accordingly, $K$ linear regression models can be established with $Y_i^{(k)}=X_i^\top\beta^{(k)}+\varepsilon_i^{(k)}$. It is notable that, both the primary task and auxiliary tasks share the same covariate vector $X_i$. Define $\bols^{(k)}\in \mR^p$ with $0\leq k \leq K$ to be the OLS estimator for each task, and let $\wh B=(\bols^{(0)},\bols^{(1)}...,\bols^{(K)})\in \mR^{p\times (K+1)}$ be the estimator matrix for all tasks. In order to borrow information from the auxiliary tasks, we develop a weighted estimator for $\beta^{(0)}$ as $\wh\beta_w = \wh Bw$, where $w \in \mR^{K+1}$ is the associated weight vector. The optimal weight $\wh w$ can be computed, which leads to the feasible weighted estimator $\wh\beta_{\wh w}$. We empirically and theoretically prove that the weighted estimator can enjoy good statistical efficiency. Finally, we apply the weighted estimator to solve the product classification problem of smart vending machines.

The rest of this article is organized as follows. Section 2 presents the estimation method for the linear regression model with auxiliary learning. The theoretical properties of the proposed weighted estimator are also rigorously studied. Section 3 presents extensive numerical experiments to demonstrate the finite sample performance of the weighted estimator. We also extend the weighted estimator to generalized linear regression models and demonstrate its performance empirically in this section. Section 4 is the real application of smart vending machines. By using the weighted estimator, the classification accuracy can be largely improved. Section 5 concludes the paper with a brief discussion.

\section{THE ESTIMATION THEORY}
\subsection{The Model Setup}

Let $(X_i,Y_i^{(0)})$ be the observation collected from the $i$-th subject with $1\leq i \leq N$. Here $Y_i^{(0)}\in\mathbb R$ is the continuous response of primary interest and $X_i=(X_{i1},...,X_{ip})^{\top}\in\mathbb R^p$ is the associated $p$-dimensional feature vector with covariance matrix $E(X_iX_i^\top)=\Sigma_{xx}\in\mR^{p\times p}$.
To model their regression relationship, a standard linear regression model can be assumed as follows:
\begin{equation}
    \label{eq:ols}
    Y_i^{(0)}=X_i^\top\beta^{(0)}+\varepsilon_i^{(0)},
\end{equation}
where $\beta^{(0)}=(\beta_1^{(0)},...,\beta_p^{(0)})^{\top}\in\mathbb R^p$ is the associated regression coefficient vector and $\varepsilon_i^{(0)}\in\mR$ is the associated random noise with $E(\varepsilon_i^{(0)})=0$ and $\var(\varepsilon_i^{(0)})=\sigma_0^2$. To estimate the regression coefficient $\beta^{(0)}$, an ordinary least squares type objective function can be constructed as $\mL^{(0)}(\beta)=N^{-1}\sum_{i=1}^N (Y_i^{(0)}-X_i^\top\beta)^2$. Let $\wh\Sigma_{xx}=N^{-1}\sum_{i=1}^N X_iX_i^\top$ and $\sxy^{(0)}=N^{-1}\sum_{i=1}^N Y_i^{(0)}X_i$. Then, an ordinary least squares (OLS) estimator can be obtained as $\bols^{(0)}=\arg\min_{\beta\in\mR^p} \mL^{(0)}(\beta)=\wh\Sigma_{xx}^{-1}\sxy^{(0)}$.
Standard asymptotic theory suggests that $E(\bols^{(0)})=\beta^{(0)}$ and $N^{-1/2}\wh\Sigma_{xx}^{1/2}(\bols^{(0)}-\beta^{(0)})\to_d \mathcal N(0,\sigma_0^2 I_p)$, where $I_p$ stands for a $p$-dimensional identity matrix  \citep{van2000asymptotic, murray2012asymptotic}. Furthermore, the finite sample performance of $\wh\beta_{\ols}^{(0)}$ is also widely demonstrated in the past literature. Specifically, it can be verified that the estimation accuracy of $\wh\beta_{\ols}^{(0)}$ measured by $\Vert\wh\beta_{\ols}-\beta^{(0)}\Vert_2$ is of the order $O_p(\sqrt{p/N})$ \citep{wainwright2019high}.

Assume we can also observe $K$ continuous variables as $Y_i^{(k)}\in\mR$ with $1\leq k\leq K$ for the $i$-th subject. These additional variables might be closely correlated to the primary response $Y_i^{(0)}$. Then we are able to adopt the idea of auxiliary learning (AL) to improve the estimation accuracy of $\beta^{(0)}$. Auxiliary learning is a machine learning technique, where a model is trained on a primary task, along with one or more auxiliary tasks that are related but not the main focus of interest. The objective of these auxiliary tasks is to provide additional guidance or regularization to the model, thus improving the estimation performance on the primary task. In auxiliary learning, $K$ is the total number of auxiliary tasks, and $Y_i^{(k)}$ with $1\leq k \leq K$ are referred to as the auxiliary responses. Assume for each auxiliary response $Y_i^{(k)}$, there exists a linear regression model as
\begin{equation}
    \label{eq:basic-regression}
    Y_i^{(k)}=X_i^\top\beta^{(k)}+\varepsilon_i^{(k)},
\end{equation}
where $\beta^{(k)}\in\mathbb R^p$ is the associated regression coefficient and $\varepsilon_i^{(k)}\in\mathbb R$ is the random noise of the $i$-th subject with $E(\varepsilon_i^{(k)})=0$ and $\var(\varepsilon_i^{(k)})=\sigma_k^2$. Let $\varepsilon_i=(\varepsilon_i^{(0)},...,\varepsilon_i^{(K)})^{\top}\in\mR^{K+1}$ be the collection of random noises for the $i$-th subject from all tasks. Assume $E(\varepsilon_i)=\bm 0_{K+1}$ and $\cov(\varepsilon_i)=\Sigma_\varepsilon$ for some covariance matrix $\Sigma_\varepsilon=(\sigma_{k_1,k_2})\in\mR^{(K+1)\times (K+1)}$. We next denote $B=(\beta^{(0)},...,\beta^{(K)})\in\mR^{p\times (K+1)}$ to collect the true regression coefficients for both primary task and auxiliary tasks.
For each $\beta^{(k)}$ with $1\leq k \leq K$, we can similarly obtain an OLS estimator as $\bols^{(k)}=\sxx^{-1}\sxy^{(k)}$, where $\sxy^{(k)}=N^{-1}\sum_{i=1}^N Y_i^{(k)}X_i$. This leads to the estimated coefficient matrix $\wh B=(\bols^{(0)},...,\bols^{(K)}) \in\mR^{p\times (K+1)}$. Meanwhile, we can obtain the estimator of the covariance matrix as $\widehat{\Sigma}_\varepsilon=(\wh\sigma_{k_1, k_2})\in\mR^{(K+1)\times (K+1)}$, where $\widehat{\sigma}_{k_1,k_2}=N^{-1}\sum_{i=1}^N (Y_i^{(k_1)}-X_i^\top\bols^{(k_1)})(Y_i^{(k_2)}-X_i^\top\bols^{(k_2)})$. Then the key question is how to exploit the information of $\wh B$ and $\widehat{\Sigma}_\varepsilon$ to improve the estimation performance of $\beta^{(0)}$.

	%%%%%%%%%%%%%%%%%%%%%%%%%%%%%%%%%%%%%%%%%%%%%%%%%%%%%%%%%%%%%%%%%%%%%%%%%%%%%%%%%%%%%%%%%%%%%%%%%%%%%%%%%%%%%%%%%%%%%%%%%%%%

	\subsection{A Weighted Estimator}

	To take advantage of the information provided by the auxiliary tasks, we study here a weighted estimator. Specifically, we assume the true coefficient matrix $B$ is of rank $d<K+1$, which potentially guarantees that $\beta^{(0)}$ can be well approximated by the linear combination of $\beta^{(k)}$s.
More specifically, consider a weight vector $w=(w_0,...,w_K)^{\top}\in\mathbb R^{K+1}$. Then let $\wh\beta_w = \wh Bw$ be a weighed estimator. Since $\wh B$ contains the OLS estimators for $B$, we should have $E(\wh B)=B$ and thus $E(\wh\beta_w) = Bw$. Recall that the ultimate goal here is to estimate the primary parameter $\beta^{(0)}$.
Consequently, our goal here is to find a weight $w$ satisfying the constraint $E(\wh\beta_w) = Bw = \beta^{(0)}$.

Let $\mathbb W=\{w\in\mathbb R^{K+1}: Bw=\beta^{(0)}\}$ be a feasible solution set. We then discuss the property of $\mathbb W$. First, we know immediately that $\mathbb W$ is not empty, since the vector $e_1=(1,0,\cdots, 0)^\top\in\mathbb R^{K+1}$ is contained in $\mathbb W$. Second, for any $w\in\mathbb W$, we should have $B(w-e_1)=\bm 0_p$. Note that $B$ is of rank $d<K+1$. Thus $\mathbb W$ should have infinite choices of $w$. Then how to find an appropriate choice of $w$ from $\mathbb W$ is an important issue. To this end, we conduct the following two steps. In the first step, we provide an analytic form of $\mathbb W$; while in the second step, under this form of $\mathbb W$, we find an optimal choice of $w$.

{\sc Step 1.} Consider the eigenvalue decomposition for $B^\top B$ as $B^\top B=\sum_{k=1}^{K+1} \lambda_k\nu_k\nu_k^\top$, where $\lambda_k$ is the $k$-th largest eigenvalue and $\nu_k=(\nu_{k1},...,\nu_{kK+1})^{\top}\\\in\mR^{K+1}$ is the associated eigenvector. Write $j_{\max}^{(k)} = \arg\max_j\{\vert\nu_{kj}\vert: 1\leq j\leq K+1\}$. By requiring $\nu_{kj_{\max}^{(k)}}>0$, both the $\lambda_k$s and $\nu_k$s can be uniquely identified.
Recall that $\rank(B)=d$ and thus $\rank(B^\top B)=d$.  We then have $\lambda_k=0$ for any $d<k\leq K+1$.
Define $\Theta=(\nu_{d+1},...,\nu_{K+1})\in\mR^{(K+1)\times (K-d+1)}$. By simple calculations, we have $\Theta^\top\Theta=I_{K+1-d}$ and $\Theta^\top B^\top B\Theta=\bm O_{K-d+1}$. Then with the help of $\Theta$, we can obtain an analytical form of  $\mathbb W$ as $\mathbb W=\{e_1+\Theta u: u\in\mR^{K-d+1}\}$. The detailed derivations are given in Appendix A.1 in the supplementary materials.
% {\color{magenta}Furthermore, for any arbitrary $u=(u_j)\in\mR^{K-d+1}$, we should have $B(e_1+\Theta u)=\beta^{(0)}$. [why??] Subsequently, we have $\mathbb W=\{e_1+\Theta u: u\in\mR^{K-d+1}\}$ [???]}

{\sc Step 2.} We next consider how to select an optimal weight $w\in\mathbb W$. Note for any $w\in\mathbb W$, we have $E(\wh\beta_w)=E(\wh Bw)=\beta^{(0)}$ and the covariance $\cov(\wh\beta_{w})=N^{-1}w^\top\Sigma_\varepsilon wE(\wh\Sigma_{xx}^{-1})$. This suggests that $\wh\beta_w$s share the same mean $\beta^{(0)}$ but with different covariance. Thus we are motivated to select an optimal weight $w^*\in\mathbb W$ that minimizes the covariance, thereby achieving the highest statistical efficiency. In other words, $w^*$ should satisfy the condition that $\cov(\wh\beta_{w})-\cov(\wh\beta_{w^*})$ is semi-positive definite for all $w\in\mathbb W$. Further define $\mathcal P(w)=w^\top \Sigma_{\varepsilon} w\in\mR$, which leads to $\cov(\wh\beta_w)=N^{-1}\mathcal P(w) E(\wh\Sigma_{xx}^{-1})$.
Note that the choice for weight $w$ only affects the term $\mathcal P(w)$ and has nothing to do with the other terms in $\cov(\wh\beta_w)$.
Therefore, the optimal choice about $w^*$ is given by $w^*=\arg\min_{w\in\mathbb W} \mathcal P(w)$. By simple calculations, we can derive
\begin{equation}
\label{eq:optimal-weight}
    w^*=e_1-\Theta(\Theta^\top \Sigma_{\varepsilon}\Theta)^{-1}\Theta^\top\Sigma_{\varepsilon} e_1.
\end{equation}
More details about the derivation of $w^*$ can be found in Appendix A.2 in the supplementary materials.

{\bf Remark 1.} By \eqref{eq:optimal-weight}, we find that the optimal weight $w^*$ is affected by the covariance $\Sigma_{\varepsilon}$. In particular, the noise level of each task (i.e., $\sigma_{k}^2$) plays an important role in determining $w^*$. To see this, consider for example a special case with $\Sigma_\varepsilon=\diag(\sigma_k^2)\in\mR^{(K+1)\times (K+1)}$ and $\beta^{(k)}=\beta^{(0)}$ for $1\leq k \leq K$. For this case, we have $w_k^*=\Big\{\sigma_k^2(\sum_{k'=0}^K \sigma_{k'}^{-2})\Big\}^{-1}$. For an arbitrary task $k$, we should have $w_k^*\to 0$ if $\sigma_k^2\to \infty$. This suggests that the optimal weights assigned to tasks with high noise levels should be small. In contrast, we should have $w_k^*\to 1$ if $\sigma_k^2\to 0$. This suggests that the optimal weights assigned to tasks with low noise levels should be large.

{\bf Remark 2.} We also find that the optimal weight $w^*$ is affected by the low-rank structure of $B$ through $\Theta$.
To see this, consider for example the case where $B$ is of full rank structure with $d=\rank(B)=K+1$. In this case, $\Theta$ becomes infeasible. Recall that for any $w\in\mathbb W$, we must have $Bw=\beta^{(0)}$ and $Be_1=\beta^{(0)}$. It follows that $B(w-e_1)=\bm 0_p$. Since $B$ is of full rank, we must have $w-e_1=\bm 0_{K+1}$. This suggests that $e_1$ is the only vector contained in the feasible set $\mathbb W$. Therefore, we must have $w^*=e_1$. In other words, the weights for all auxiliary tasks should be zero, i.e., $w_k^*=0$ for $1\leq k\leq K$. Therefore all the auxiliary tasks are useless and should be excluded.

By plugging in the optimal weight $w^*$ into $\wh\beta_w$, we obtain the estimator $\wh\beta_{w^*}=\wh Bw^*$, which is referred to as the ``ORACLE" estimator. However, the ORACLE estimator cannot be used immediately, since the optimal weight $w^*$ involves two unknown parameters (i.e., $\Theta$ and $\Sigma_\varepsilon$).
To find a practically feasible estimator, both $\Theta$ and $\Sigma_\varepsilon$ need to be estimated. In Section 2.1, we have already obtained an estimator for the covariance matrix as $\widehat{\Sigma}_\varepsilon$. Therefore the key problem here is to find an estimator for $\Theta$. To this end, we seek the help from the estimated coefficient matrix $\wh B$. Specifically, we can conduct the eigenvalue decomposition on $\wh B^\top \wh B$. Let $\wh\lambda_k$ be the $k$-th largest eigenvalue and $\wh\nu_k$ be the associated eigenvector. Therefore, $w^*$ can be estimated as
$\wh w^*=e_1-\wh\Theta (\wh\Theta^\top \wh\Sigma_\varepsilon\wh\Theta)^{-1} \wh\Theta^\top \wh\Sigma_\varepsilon e_1$, where $\wh\Theta=(\wh\nu_{d+1},...,\wh\nu_{K+1}) \in \mR^{(K+1)\times (K-d+1)}$.
We then plug $\wh w^*$ into $\wh\beta_{w^*}$ and this leads to the final weighted estimator $\wh\beta_{\wh w^*}$, which is referred to as the ``FEASIBLE" estimator.

\subsection{Convergence Rate}
In this section, we study the theoretical properties of the feasible estimator $\wh\beta_{\wh w^*}$. To facilitate the theoretical analysis, we first introduce some useful notations. For any arbitrary matrix $A$, define the operator norm as $\Vert A\Vert_{\op}=\lambda_{\max}^{1/2}(A^\top A)$, where $\lambda_{\max}(M)$ stands for the maximum absolute eigenvalue of an arbitrary symmetric matrix $M$.
Furthermore, denote the Frobenius norm for $A$ as $\Vert A\Vert_F=\tr^{1/2}(A^\top A)$. Moreover, for any arbitrary vector $u$, define its $\ell_2$-norm as $\Vert u\Vert_2=(u^\top u)^{1/2}$.
Next, for a random variable $U\in\mathbb R$, define its sub-Gaussian norm as $\Vert U\Vert_{\psi_2}=\inf\{t>0: E\{\exp(U^2/t^2)\}\leq 2\}$. For a random vector $V\in\mathbb R^p$, define its sub-Gaussian norm as $\Vert V\Vert_{\psi_2}=\sup_{\Vert u\Vert = 1}\Vert u^\top V \Vert_{\psi_2}$. We then say a random vector $V$ is a sub-Gaussian random variable if $\Vert V\Vert_{\psi_2}<\infty$. Please see \citet{vershynin2018high} and \citet{wainwright2019high} for more discussions about the theoretical properties of sub-Gaussian variables.
To study the non-asymptotic behavior of the proposed FEASIBLE estimator, some commonly used technical conditions are needed and given as follows.
\begin{itemize}
    \item[(C1)] {\sc (Sub-Gaussian Condition)} Assume $V_i=(X_i^\top, \varepsilon_i^\top)^\top \in\mathbb R^{p+K+1}$ is an independent and identically distributed sub-Gaussian random vector with $\Vert V_i\Vert_{\psi_2}=C_{\sub}<\infty$.

    \item[(C2)] {\sc (Variance Condition)} Assume there are two constants $\tau_{\min}$ and $\tau_{\max}$ such that (1) $0<\tau_{\min}\leq \lambda_{\min}(\Sigma_{xx})\leq \lambda_{\max}(\Sigma_{xx})\leq \tau_{\max}<\infty$ as $p\to\infty$; and (2) $0< \tau_{\min}\leq \lmin(\Sigma_{\varepsilon}) \leq \lmax(\Sigma_{\varepsilon})\leq \tau_{\max}<\infty$ as  $K\to\infty$.

    \item[(C3)] {\sc (Regression Coefficients)} Assume $\rank(B)=d<K+1$. For the same constants $0<\tau_{\min}\leq \tau_{\max}<\infty$ as in Condition (C2), we assume that $\tau_{\min}\leq\lambda_{d}(B^\top B)\leq \lambda_{1}(B^\top B)=\lambda_{\max}(B^\top B)\leq\tau_{\max}$ as $p\to\infty$. Here $\lambda_{j}(A)$ stands for the $j$-th largest eigenvalue for an arbitrary symmetric matrix $A$.

    \item[(C4)] {\sc (Divergence Rate)} Assume that (1) $K\to\infty$, (2) $p\to\infty$, (3) $K/p\to 0$, and (4) $p/N\to 0$ as $N\to\infty$.

\end{itemize}

Condition (C1) is a standard sub-Gaussian assumption commonly used in high-dimensional researches \citep{vershynin2018high, wainwright2019high, zhu2022feature, tian2023transfer}.
Condition (C2) imposes a fixed lower bound for $\lambda_{\min}(\Sigma_{xx})$ so that different features are not allowed to be linearly dependent with each other severely. In the meanwhile, Condition (C2) imposes a fixed upper bound on $\lambda_{\max}(\Sigma_{xx})$ so that the variance of $X_i$ is bounded.
The second condition in (C2) is used to prevent the random noises of different tasks from being seriously correlated. Conditions of similar type have been widely used in the past literature \citep{fan2008high, zhu2020multivariate, fan2024latent}.
Condition (C3) regularizes the low-rank structure for the regression coefficient $B$. On one side, it provides an upper bound for $\lambda_{\max}(B^\top B)$, so that the information contained in $X_i\beta^{(k)}$s would not explode to infinity. On the other side, the lower bound of $\lambda_{d}(B^\top B)$ enforces that different $\beta^{(k)}$s should not be linearly dependent with each other severely.
Condition (C4) specifies the divergence rates of $K$, $p$, and $N$.
Under the above conditions, we have the following theorem.
\begin{theorem}
\label{thm-1}
    Assume the technical conditions (C1)---(C4) hold. Then there exist constants $C_1$ and $C_2$, such that the following inequality holds with probability at least $1-C_2\exp(-K)$, as long as $N>N_0$ for some sufficiently large constant $N_0$
    \begin{equation}
        \label{eq:thm-1}
        \Big\Vert\wh\beta_{\wh w^*}-\wh\beta_{w^*}\Big\Vert_2\leq C_1\left\{\sqrt{\frac{K+d}{N}}+\frac{p}{N}\right\}.
    \end{equation}
\end{theorem}
\noindent
The detailed proof of Theorem \ref{thm-1} is given in Appendix B in the supplementary materials.
Theorem \ref{thm-1} provides an upper bound for the discrepancy between the FEASIBLE estimator $\wh\beta_{\wh w^*}$ and the ORACLE estimator $\wh\beta_{w^*}$. From Theorem \ref{thm-1}, we know immediately that $\Vert\wh\beta_{\wh w^*}-\wh\beta_{w^*}\Vert_2=O_p\big\{\sqrt{(K+d)/N}+p/N\big\}$ as $N\to\infty$. This finding suggests that the discrepancy between $\wh\beta_{\wh w^*}$ and $\wh\beta_{w^*}$ goes to zero as $N\to\infty$.

Based on Theorem \ref{thm-1}, we further study the discrepancy between the FEASIBLE estimator $\wh\beta_{\wh w^*}$ and the true parameter $\beta^{(0)}$. To this end, we first study the theoretical behaviour of the ORACLE estimator $\wh\beta_{w^*}$. Note that $\wh\beta_{w^*}$ is an unbiased estimator for $\beta^{(0)}$, i.e., $E(\wh\beta_{w^*})=\beta^{(0)}$. Moreover, since $\wh\beta_{w^*}$ is a linear combination of several OLS estimators, it can be verified that $\sqrt{N}\alpha^\top (\wh\beta_{w^*}-\beta^{(0)})\to_d \mathcal N\Big(0, \mathcal P(w^*)\alpha^\top\Sigma_{xx}^{-1}\alpha\Big)$, where $\alpha=(\alpha_1,...,\alpha_p)^{\top}\in\mR^p$ is an arbitrary $p$-dimensional, non-zero, and fixed vector. This suggests that, the difference between $\wh\beta_{w^*}$ and $\beta^{(0)}$ is of the order $O_p(\sqrt{p/N})$, in the sense that $\alpha^\top(\wh\beta_{w^*}-\beta^{(0)})=O_p(\sqrt{p/N})$. On the other side, Theorem \ref{thm-1} suggests that $\alpha^\top(\wh\beta_{\wh w^*}-\wh\beta_{w^*})=O_p\{\sqrt{(K+d)/N} + p/N\}$ as $N\to\infty$. It is a small order term as compared with $\alpha^\top(\wh\beta_{w^*}-\beta^{(0)})=O_p(\sqrt{p/N})$ since $d< K+1$, $K/p\to 0$, and $p/N\to 0$ as $N\to\infty$. Therefore, the difference between $\wh\beta_{\wh w^*}$ and $\wh\beta_{w^*}$ can be asymptotically ignorable. This suggests that $\alpha^\top(\wh\beta_{\wh w^*}-\beta^{(0)})=\alpha^\top(\wh\beta_{\wh w^*}-\wh\beta_{w^*}) + \alpha^\top(\wh\beta_{w^*}-\beta^{(0)})$ should have the same asymptotic distribution as $\alpha^\top(\wh\beta_{w^*}-\beta^{(0)})$. Therefore, we have the following corollary. 

\begin{corollary}
\label{coro1}
Assume the technical conditions (C1)---(C4) hold. Then as $N\to\infty$, we have $\sqrt{N}\alpha^\top (\wh\beta_{\wh w^*}-\beta^{(0)})\to_d \mathcal N\Big(0, \mathcal P(w^*)\alpha^\top\Sigma_{xx}^{-1}\alpha\Big)$.
\end{corollary}

\section{SIMULATION STUDY}

\subsection{Simulation Setup and Performance Measure}

To demonstrate the finite sample performance of the weighted estimators, we conduct a number of simulation experiments in this section. We start with the primary linear regression model \eqref{eq:ols}. The feature vector $X_i$ with $1\leq i \leq N$ is generated from a multivariate Gaussian distribution with mean $\bm{0}_p$ and covariance matrix $\Sigma_{xx} = (\sigma_{j_1,j_2})\in\mR^{p\times p}$, where $\sigma_{j_1,j_2} = 0.5^{|j_1 - j_2|}$ for $1\leq j_1, j_2 \leq p$ \citep{tibshirani1996regression, fan2008sure}. Recall that $B \in \mathbb R ^{p\times (K+1)}$ is the true coefficient matrix with rank $d$.
We generate $B$ as $B=\Phi_1\Phi_2/\Vert\Phi_1\Phi_2\Vert_{\op}$ with $\Phi_1\in \mathbb R ^{p\times d}$ and $\Phi_2\in \mathbb R ^{d\times (K+1)}$. Here $\Phi_1=(\phi_{ij}^{(1)})$ with $\phi_{ij}^{(1)}$ independently generated from $\mathcal N(0,1/\sqrt p)$;
$\Phi_2=(\phi_{ij}^{(2)})$, where $\phi_{ij}^{(2)}=0$ if $i\leq j\leq i+K-d+1$ and $\phi_{ij}^{(2)}=1$ otherwise.
It follows that $\rank(B)=d$ and $\Vert B\Vert_{\op}=1$.
The random noise vector $\varepsilon_i=(\varepsilon_i^{(0)},...,\varepsilon_i^{(K)})^{\top}\in\mR^{K+1}$ is generated independently from a multivariate normal distribution with mean $\bm 0_{K+1}$ and covariance matrix $\Sigma_\varepsilon\in\mR^{(K+1)\times (K+1)}$. Here $\Sigma_\varepsilon=\wt\Sigma_\varepsilon\wt\Sigma_\varepsilon^{\top}/\Vert\wt\Sigma_\varepsilon\Vert_{\op}^2$ and each component of $\wt\Sigma_\varepsilon\in\mR^{(K+1)\times (K+1)}$ is independently generated from $\mathcal N(0,1/\sqrt{K+1})$. It follows that $\Vert\Sigma_\varepsilon\Vert_{\op}=1$. 
Based on the feature vector $X_i$, the coefficient matrix $B$, and the random noise $\varepsilon_i$, the responses $Y_i^{(k)}$ with $0\leq k \leq K$ can be generated accordingly.

We randomly replicate the data generation process for a total of $M=500$ times. For each generated dataset, we compute three estimators for the primary coefficient $\beta^{(0)}$. The first one is the OLS estimator by using $Y_i^{(0)}$ and $X_i$ only. We denote this estimator by ``OLS". The second one is the ORACLE estimator $\wh\beta_{w^*}$ using the optimal weight $w^*$. The last one is the FEASIBLE estimator $\wh\beta_{\wh w^*}$ using the estimated weight $\wh w^*$. For each estimator (i.e., OLS, ORACLE, and FEASIBLE), we define $\widehat{\beta}_{(m)}^{(0)} \in \mathbb R^{p}$ as the estimator for $\beta^{(0)}$ in the $m$th replication ($1\leq m \leq M$). Then, to evaluate the estimation performance of each estimator, we use the mean squared error (MSE), namely, $M^{-1}\sum_{m=1}^{M}\|\widehat{\beta}_{(m)}^{(0)}-\beta^{(0)}\|_2^2$, where $\|\cdot \|_2$ is the $\ell_2$-norm.

\subsection{Example 1: Varying $N$ and $p$}

\begin{figure}[h]
    \centering
    \subfloat[$p=\sqrt{N},K=10,d=5$]{
        \includegraphics[width=.48\textwidth]{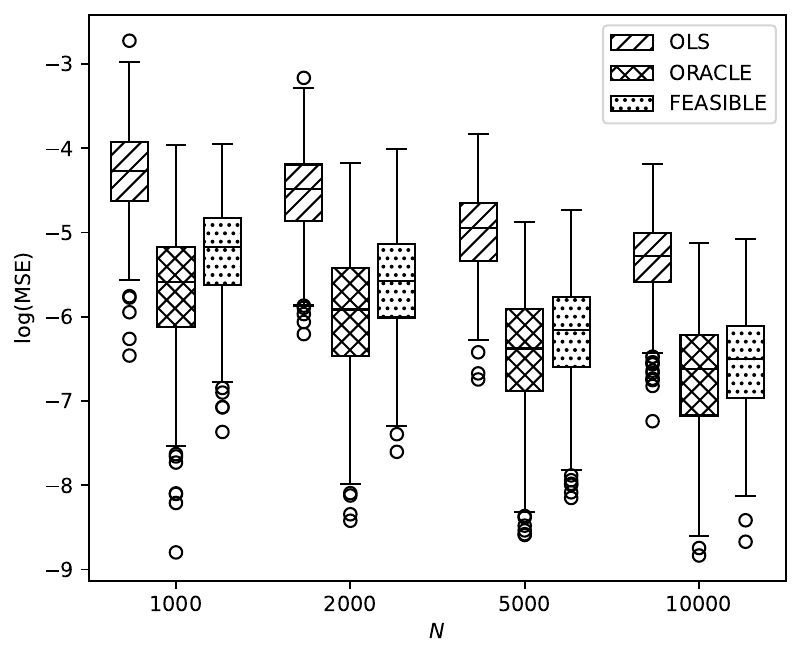}
    }
    \subfloat[$p=\sqrt{N},K=20,d=10$]{
        \includegraphics[width=.48\textwidth]{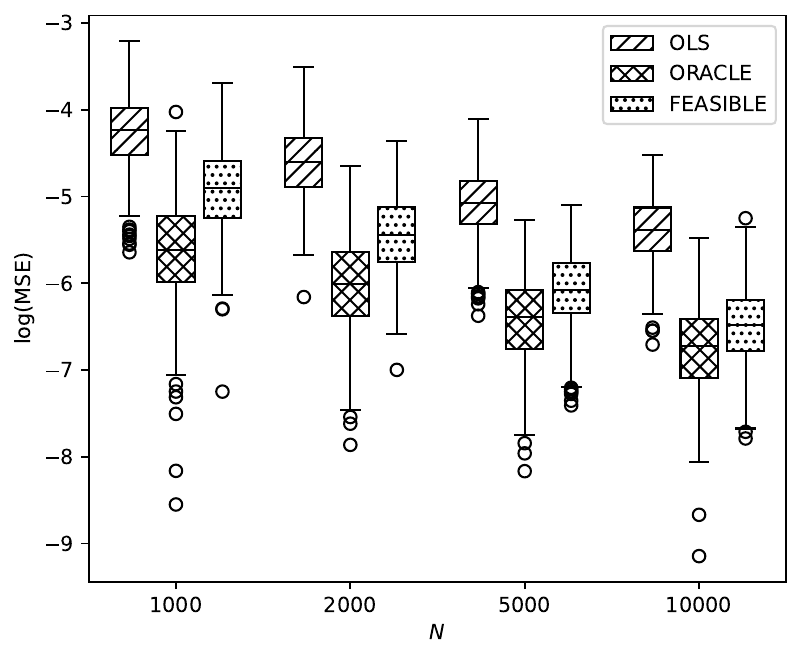}
    }
    \caption{Boxplots of log-transformed MSE values of the OLS, ORACLE, and FEASIBLE estimators with varying $N$ and $p$.}
    \label{fig:exp1_unfix_np}
\end{figure}

We first study the estimation performance of different estimators with varying sample size $N$ and feature dimension $p$. Specifically, we assume the total sample size $N\in \{1000,2000,5000,\\ 10000\}$. Let the feature dimension $p = \sqrt{N}$, so that the feature dimension $p$ diverges as the sample size $N$ increases. Then for each $(N,p)$-pair, fix the number of auxiliary tasks $K\in \{10,20\}$ and let $d=\rank(B)=K/2$ respectively. The detailed simulation results are shown in Figure \ref{fig:exp1_unfix_np}, from which we can draw the following conclusions. First, when the sample size (i.e., $N$) increases, the MSE values of the three estimators steadily decrease. This suggests the consistency of these estimators. Second, regardless of the sample size, the MSE values of the ORACLE and FEASIBLE estimators are consistently smaller than those of the OLS estimator. This finding demonstrates the advantage of the weighted estimators over the OLS estimator with the help of auxiliary tasks. Third, when $N$ is relatively large, the FEASIBLE estimator can achieve almost the same estimation performance as the ORACLE estimator, which is consistent with Theorem \ref{thm-1}.

\subsection{Example 2: Varying $K$ with Fixed $d$}

\begin{figure}[h]
    \centering
    \subfloat[$N=2000,p=\sqrt{N},d=5$]{
        \includegraphics[width=.48\textwidth]{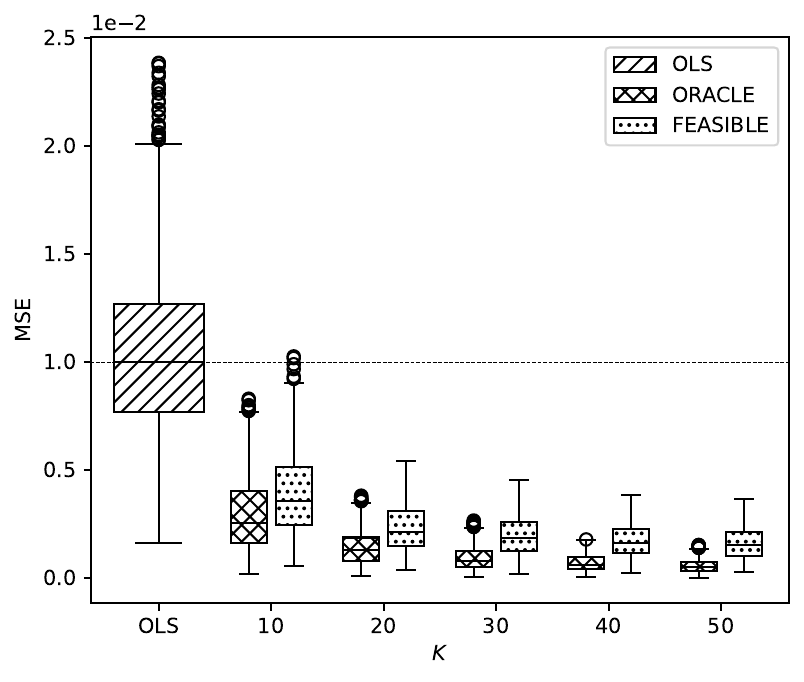}
    }
    \subfloat[$N=10000,p=\sqrt{N},d=5$]{
        \includegraphics[width=.48\textwidth]{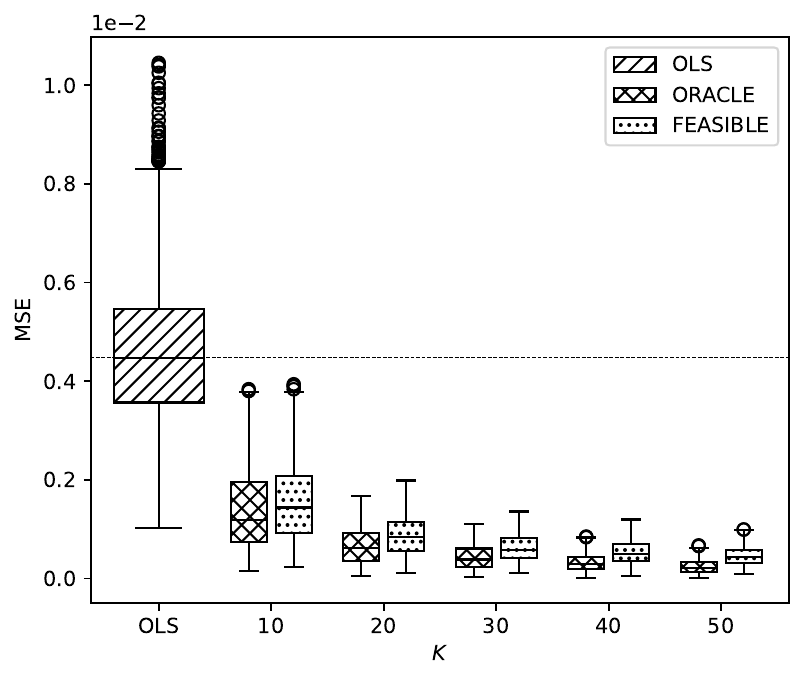}
    }
    \caption{Boxplots of MSE values of the OLS, ORACLE, and FEASIBLE estimators with fixed $d$ but varying $K$.}
    \label{fig:exp1_unfix_d}
\end{figure}

We next investigate how the number of auxiliary tasks affects the estimation performance of the weighted estimators. In this scenario, we assume that $d = \text{rank}(B)$ remains fixed. In other words, as the number of auxiliary tasks $K$ increases, the newly added coefficient vectors remain within the original column space of $B$. In this case, the auxiliary tasks are highly correlated with each other. It would lead to a more flexible feasible solution set $\mathbb W$, which further yields a smaller $\mathcal P(w^*)=\min_{w\in\mathbb W} w^\top\Sigma_\varepsilon w$. Then by our Corollary \ref{coro1}, we should expect both ORACLE and FEASIBLE estimators have lower MSE values as $K$ increases. To illustrate this idea, assume $N\in \{2000,10000\}$ and $p = \sqrt{N}$. Let $K$ vary from 10 to 50 with a step size of 10. For each $K$, we fix $d=5$. The corresponding simulation results are shown in Figure \ref{fig:exp1_unfix_d}. We can find that, with a fixed $d$, the MSE values of both the ORACLE and FEASIBLE estimators decrease as $K$ grows.
This finding also verifies that under a fixed $d$, including more auxiliary tasks should improve the estimation performance.

\subsection{Example 3: Varying $d$ with Fixed $K$}

\begin{figure}[h]
    \centering
    \subfloat[$N=2000,p=\sqrt{N},K=10$]{
        \includegraphics[width=.48\textwidth]{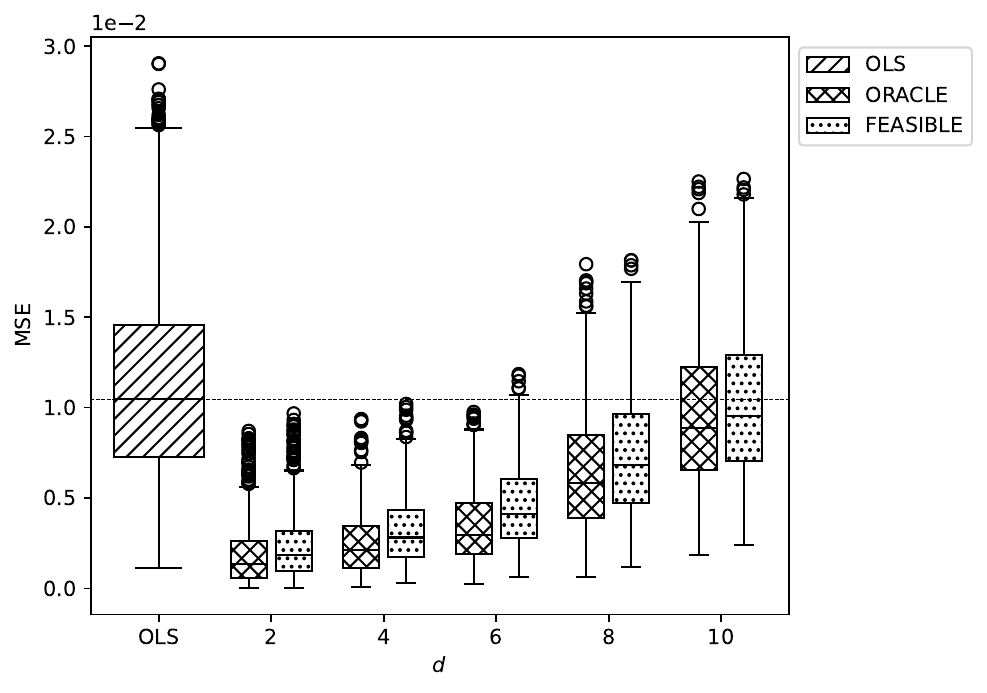}
    }
    \subfloat[$N=10000,p=\sqrt{N},K=10$]{
        \includegraphics[width=.48\textwidth]{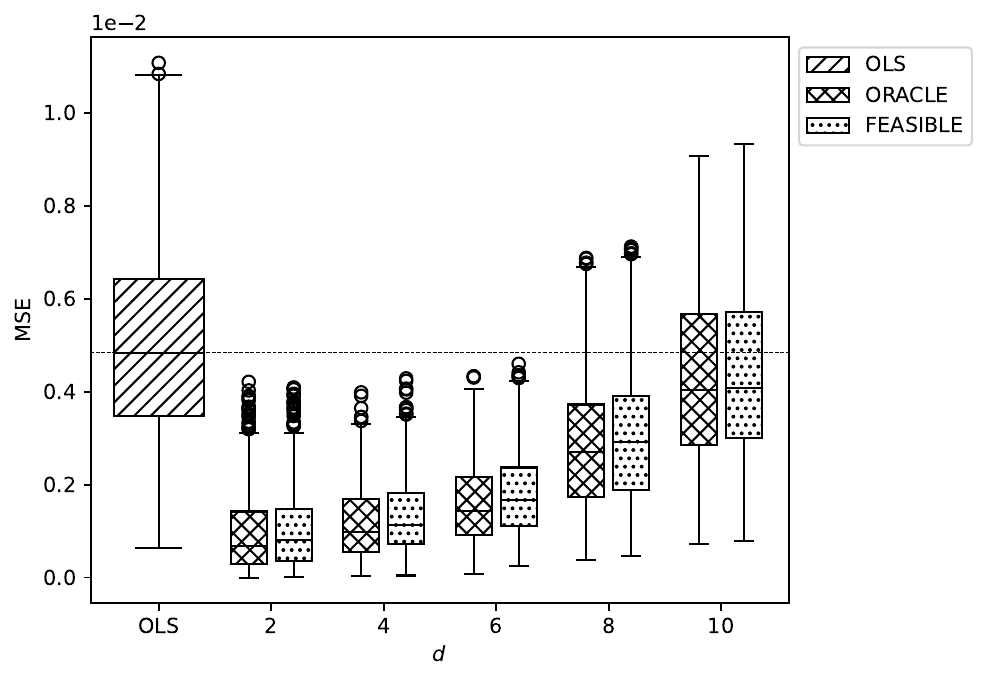}
    }
    \caption{Boxplots of MSE values of the OLS, ORACLE, and FEASIBLE estimators with fixed $K$ but varying $d$.}
    \label{fig:exp1_unfix_k}
\end{figure}

In this example, we study the effect of the rank $d$ under a fixed number of auxiliary tasks $K$. With a fixed $K$, a larger rank $d$ suggests $\beta^{(k)}$s are less likely to be linearly dependent with each other, and thus might introduce extra noise and harm the estimation performance.
To verify this idea, we assume $N\in \{2000,10000\}$, $p = \sqrt{N}$, and fix $K=10$. Then we investigate the estimation performance of the weighted estimators with varying $d$. Specifically, let $d$ vary from 2 to 10 with a step size of 2. The corresponding results are shown in Figure \ref{fig:exp1_unfix_k}.
By Figure \ref{fig:exp1_unfix_k}, it can be verified that with an increasing $d$, the MSE values for both ORACLE and FEASIBLE estimators are increasing steadily, which is consistent with our expectation.

\subsection{Example 4: Estimating with Low-Quality Auxiliary Tasks}

In this example, we focus on how the quality of auxiliary tasks affects the estimation performance. Assume we already have a total of $K$ auxiliary tasks, whose coefficient matrix is denoted by $B_1$ and $\rank(B_1)=d$. Assume all these auxiliary tasks are beneficial to the primary task, with each having an optimal weight $w_k^*\neq 0$ for every $1\leq k\leq K$. We next introduce an additional set of $K'$ auxiliary tasks, whose coefficient matrix is denoted by $B_2=(\beta^{(K+1)},...,\beta^{(K+K')})\in \mR^{p\times K'}$. Further assume these newly introduced auxiliary tasks are of low-quality relative to the primary task, i.e., their optimal weights are given by $w^*_{k}=0$ for $K<k\leq K+K'$. 
In summary, we have a total of $K+K'$ auxiliary tasks, where the first $K$ tasks are useful and the remaining $K'$ tasks are useless.
Next, we fix $N=10,000$, $p=\sqrt N=100$, $K=10$, $d=5$, and $K'=50$. 
We then generate the coefficient matrix as follows. Similar to the generation process in Section 3.1, we first separately generate two coefficient matrices $B_1\in\mR^{p\times (K+1)}$ with $\rank(B_1)=d$ and $B_2\in\mR^{p\times K'}$ with $\rank(B_2)=K'$. Then we define $B_2^* = \{I- B_1(B_1^\top B_1)^{-1} B_1^\top\}B_2$, so that $B_2^*$ is linearly independent of $B_1$. This leads to the final coefficient matrix $B=(B_1, B_2^*)\in\mR^{p\times (K+K'+1)}$. Subsequently, the covariate $X_i$ and random noise $\varepsilon_i$ are generated similar to Section 3.1, and the response $Y_i$ can also be generated accordingly.

\begin{figure}[h]
    \centering
    \subfloat{
        \includegraphics[width=.60\textwidth]{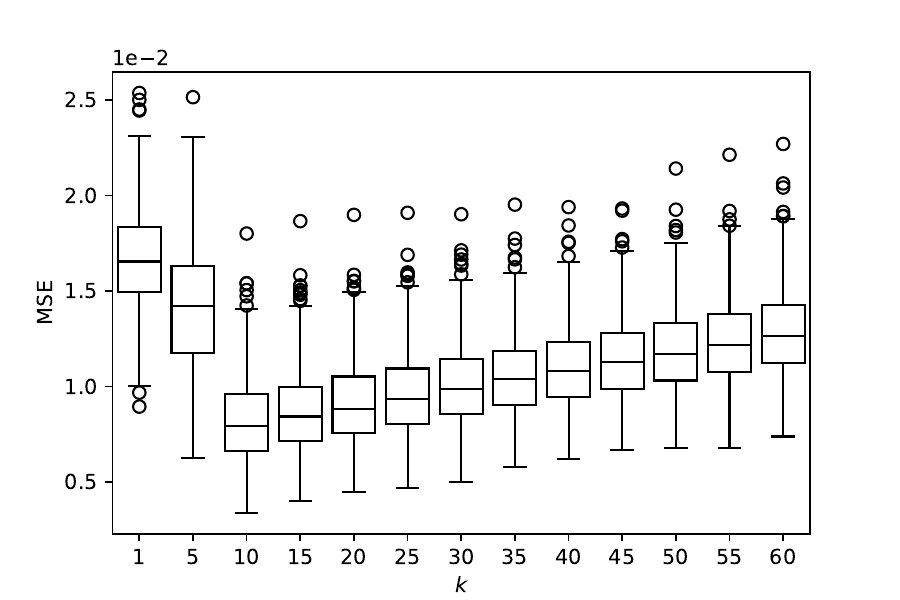}
    }
    \caption{Boxplots of MSE values of FEASIBLE estimators using different numbers of auxiliary tasks.}
    \label{fig:exp2_lowqual}
\end{figure}

To evaluate the estimation performance by incorporating different number of auxiliary tasks, we denote $\wh\beta_{\wh w^*}(k)$ as the FEASIBLE estimator using the first $k$ auxiliary tasks. For example, $\wh\beta_{\wh w^*}(10)$ is the FEASIBLE estimator only using all useful auxiliary tasks, and $\wh\beta_{\wh w^*}(k)$ with $k>10$ is the FEASIBLE estimator using both useful and useless auxiliary tasks. The MSE values are then computed and reported in Figure \ref{fig:exp2_lowqual}. From Figure \ref{fig:exp2_lowqual}, an interesting $U$-shape pattern is detected. Specifically, we find that the MSE values decrease as $k$ increases for $1\leq k\leq 10$. This is expected since the auxiliary tasks included at this stage are all useful with $w_k^*\neq 0$.
However, once $k>10$, the MSE value starts to increase steadily as $k$ increases. This is also expected since the newly included tasks at this stage are of low quality with $w_k^*=0$.
This suggests that including low-quality tasks might lead to degraded estimation efficiency.

\subsection{An Extension: Logistic Regression}

In this section, we evaluate the possibility of extending our method from a linear regression model to some other models. To illustrate this idea, we consider the popularly used logistic regression model as an example. We first follow the data generation process in Section 3.1 to generate the feature vector $X_i$ and the coefficient matrix $B$. For simplicity, we assume a diagonal covariance matrix $\Sigma_\varepsilon$ for the random noise. That is, $\Sigma_\varepsilon=\diag(\sigma_k^2:0\leq k\leq K)\in\mR^{(K+1)\times (K+1)}$. Then assume the response $Y_i^{(k)}$ with $0\leq k \leq K$ is a binary variable, which is generated from a Bernoulli distribution with the probability given as $p_{ik}=P(Y_i^{(k)}=1|X_i,\beta^{(k)})=\exp(X_i^{\top}\beta^{(k)})/\{1+\exp(X_i^{\top}\beta^{(k)})\}$. The data generation process is replicated for $M=500$ times. Then we compare the ORACLE and FEASIBLE estimators with the maximum likelihood estimator (MLE) only using the primary response $Y^{(0)}_i$ and $X_i$. For the ORACLE and FEASIBLE estimators, the weights $w^*$ and $\wh w^*$ are calculated similarly as those in the linear regression model. The only difference lies in the computation of variance terms. That is $\wh\Sigma_\varepsilon=\diag(\wh\sigma_k^2)\in\mR^{(K+1)\times (K+1)}$, where
$\widehat{\sigma}_k^2=N^{-1}\sum_{i=1}^N \widehat{p}_{ik}(1-\wh p_{ik})$. Here $\widehat{p}_{ik}$ is the estimated counterpart of $p_{ik}$ with $\beta^{(k)}$s replaced by the corresponding MLE.

\begin{figure}[h]
    \centering
    \subfloat[Varying $N$ and $p$]{
        \includegraphics[width=.34\textwidth]{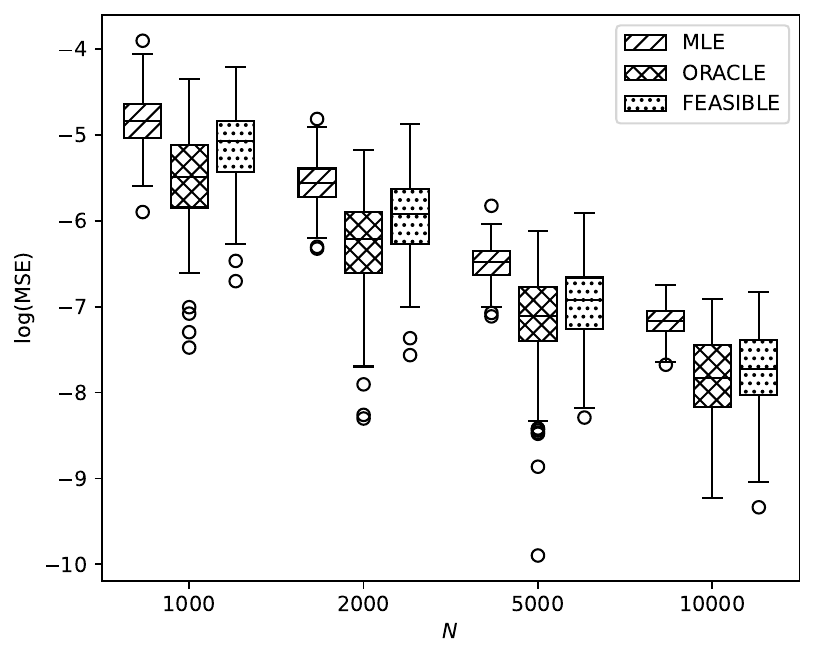}
    }
    \subfloat[Varying $K$, Fixed $d$]{
        \includegraphics[width=.32\textwidth]{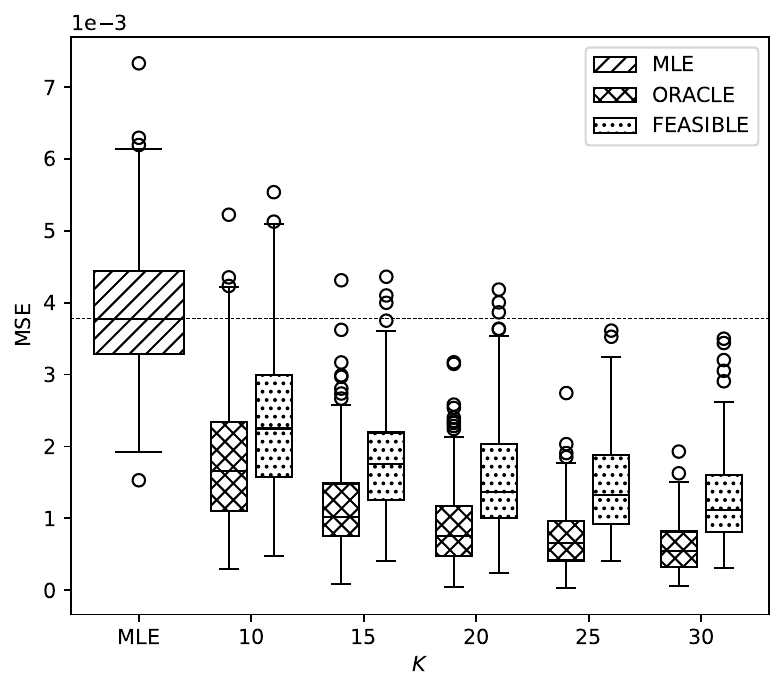}
    }
    \subfloat[Varying $d$, Fixed $K$]{
        \includegraphics[width=.32\textwidth]{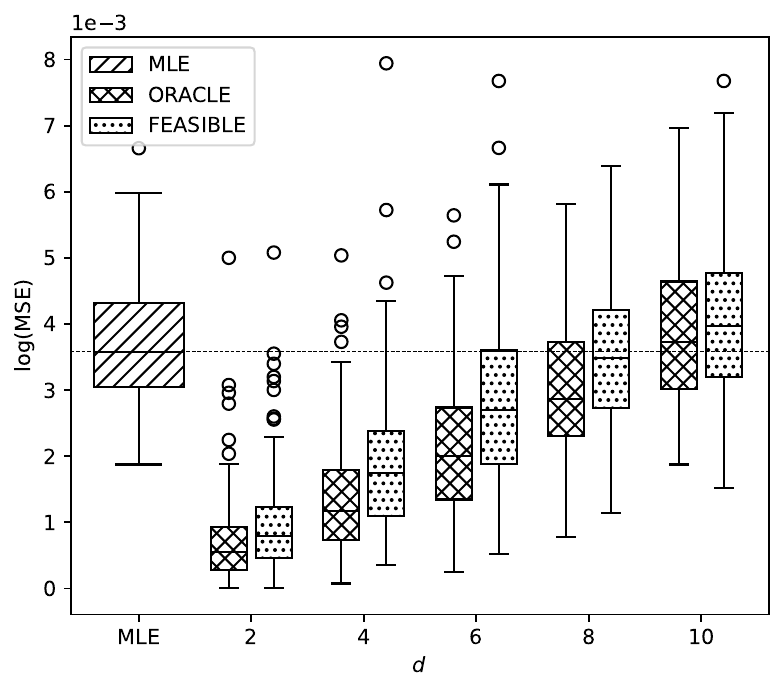}
    }
    \caption{Boxplots of MSE values of the MLE, ORACLE, and FEASIBLE estimators for the logistic regression model. For varying $N$ and $p$, we fix $K=10$ and $d=5$. For varying $K$, we fix $d=5$, $N=2000$, and $p=\sqrt{N}$. For varying $d$, we fix $K=10$, $N=2000$, and $p=\sqrt{N}$}
    \label{fig:exp2_lr_unfix_np}
\end{figure}

We consider different combinations of $(N,p,K,d)$ to explore the estimation performance of the weighted estimators. Specifically, for varying $N$ and $p$, we fix $K=10$ and $d=5$. For varying $K$ and fixed $d$, we assume $N=2000$, $p=\sqrt{N}$, and $d=5$. For varying $d$ and fixed $K$, assume $K=10$, $N=2000$, and $p=\sqrt{N}$. The detailed simulation results are shown in Figure \ref{fig:exp2_lr_unfix_np}. In general, the results under the logistic regression model are quantitatively similar to those under the linear regression model. We demonstrate that the ORACLE and FEASIBLE estimators are both consistent as the MSE values decrease with the increase of the sample size. The two weighted estimators also enjoy better statistical efficiency than the MLE as $B$ is rank-reducible.

\section{REAL DATA ANALYSIS}

\subsection{Data Description}

\begin{figure}[h]
    \centering
    \subfloat[{\centering Smart vending machine}]{
        \includegraphics[height=.3\textheight]{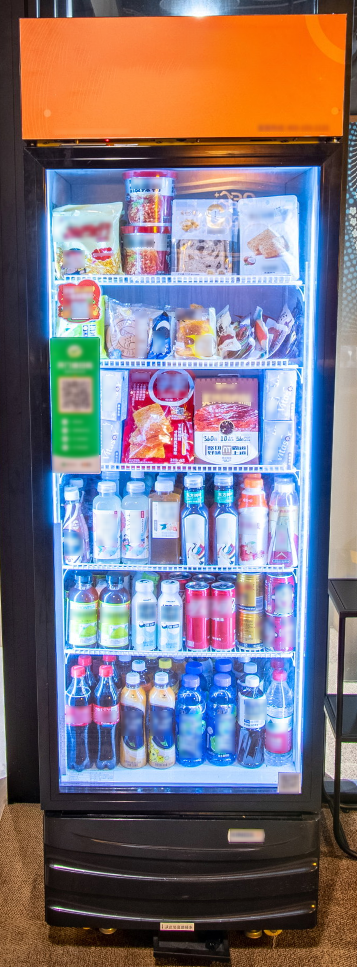}
    }
    \subfloat[45 types of bottled water]{
        \includegraphics[height=.3\textheight]{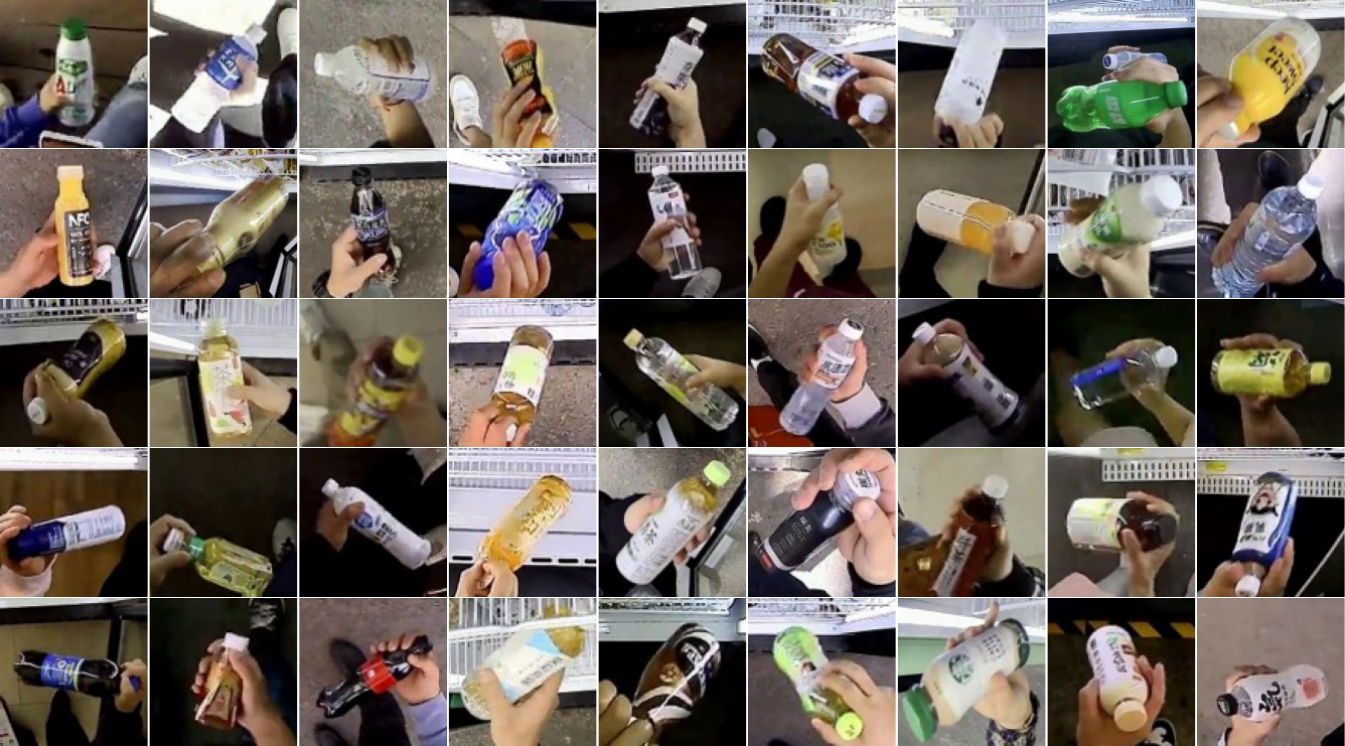}
    }
    \caption{The left panel is one smart vending machine and the right panel shows 45 types of bottled water sold at the vending machine.}
    \label{fig:vending-machine-primary-task}
\end{figure}

To demonstrate the practical performance of our method, we conduct an empirical study to identify products in smart vending machines, which are commonly encountered in our daily life. Here we focus on one particular type of smart vending machine as shown in the left panel of Figure \ref{fig:vending-machine-primary-task}. With this type of smart vending machine, users can experience a convenient shopping process characterized by ``scan to open, choose freely, and automatically settle when closed" \citep{Wang2022Design}. Each smart vending machine is equipped with a camera that records the shopping process of users. Then based on the video recordings, the smart vending machine needs to accurately identify the products purchased by the user and then establish financial settlement. Thus in this work, our primary goal is to identify the products. For illustration purposes, we choose the product type of ``bottled water", which covers a total of forty-five categories (see the right panel of Figure \ref{fig:vending-machine-primary-task}). The dataset is provided by a Chinese smart retail company, which consists of $N=6,018$ transactions for bottled water. For each transaction, a sub-image containing the purchased product is carefully sampled from the camera recordings by an appropriately designed automatic approach. This leads to $N=6,018$ photos, each of which corresponds to a unique transaction. In the subsequent analysis, we aim to identify the products based on the photo information.

\subsection{The Design of Auxiliary Tasks}

\begin{figure}[h]
    \centering
    \subfloat{
        \includegraphics[height=.25\textheight]{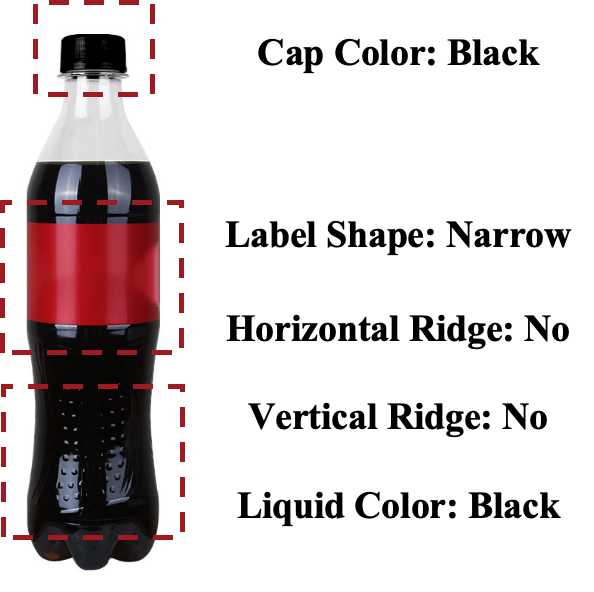}
    }
    \subfloat{
        \includegraphics[height=.25\textheight]{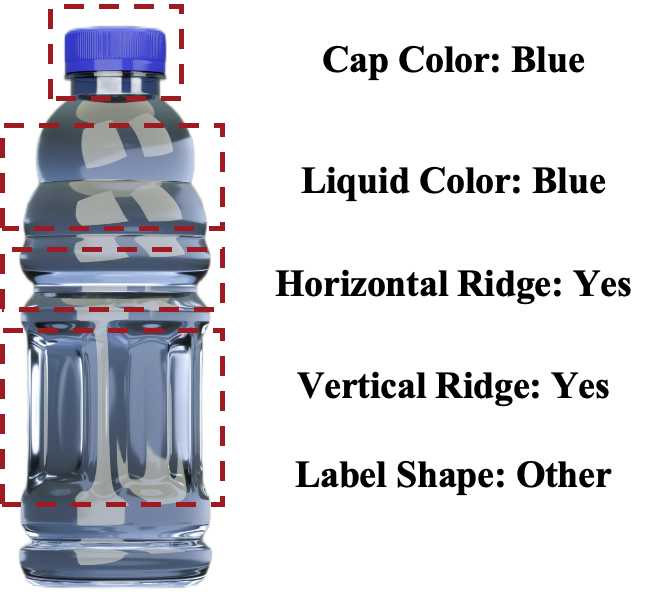}
    }
    \caption{The illustration of two specific products. 
    The product on the left panel is with ``black'' cap color, ``narrow'' label shape, and ``black'' liquid color. It does not have horizontal or vertical ridges. The product on the right panel is with ``blue'' cap color, ``blue'' liquid color, and ``other'' label shape. It has both horizontal and vertical ridges.}
    \label{fig:auxiliary-example}
\end{figure}

\begin{figure}
    \centering
    \subfloat[Cap colors with ten classes]{
        \includegraphics[width=.8\textwidth]{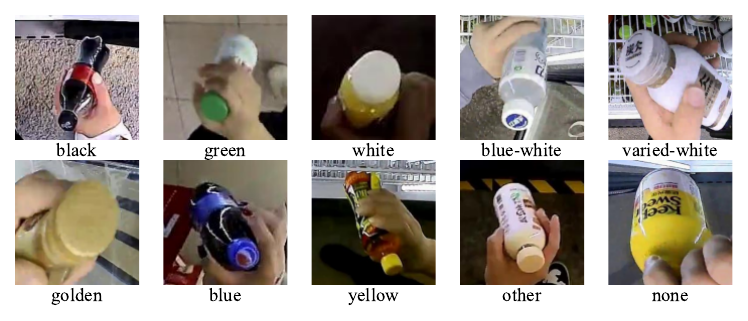}
    }

    \subfloat[Label shapes with five classes]{
        \includegraphics[width=.8\textwidth]{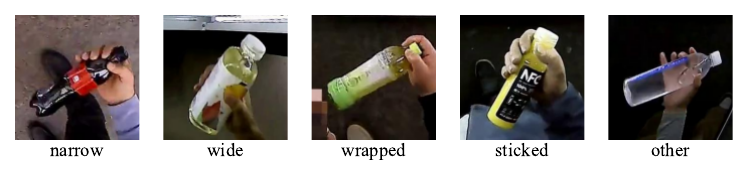}
    }

    \subfloat[Liquid colors with nine classes ]{
        \includegraphics[width=.8\textwidth]{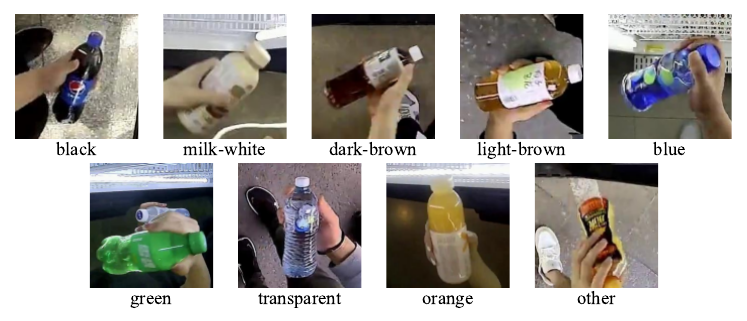}
    }

    \subfloat[Horizontal ridge]{
        \includegraphics[width=.33\textwidth]{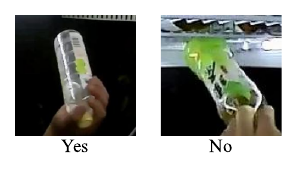}
    }
    \subfloat[Vertical ridge]{
        \includegraphics[width=.33\textwidth]{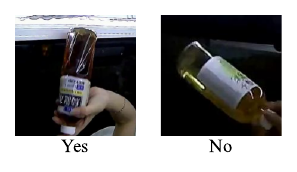}
    }
    \caption{The detailed classification categories of the five auxiliary tasks. Each item has a label in each of the five auxiliary classification tasks. }
    \label{fig:auxiliary}
\end{figure}

We treat this product classification problem as a multi-class classification task. Specifically, define the primary response as the product category labeled as $Y_i^{(0)}\in\{1,\cdots,45\}$. Based on some preliminary exploration, directly classifying $Y_i^{(0)}$ does not yield satisfactory results. To improve the classification accuracy, we seek the help of auxiliary learning. Specifically, we consider five auxiliary tasks. The first one is identifying the cap color of the bottle (e.g., white, blue, red). It is a ten-class task and thus  $Y_i^{(1)} \in\{1,\cdots,10\}$. The second one is identifying the label shape of the bottle (e.g., wide, narrow, wrapped). It is a five-class task and thus we have $Y_i^{(2)}\in\{1,\cdots,5\}$. The third one is identifying the liquid color (e.g., white, black, yellow). It is a nine-class task and thus $Y_i^{(3)}\in\{1,\cdots,9\}$.
The fourth task is to identify whether there are horizontal ridges on the plastic bottle. This is a binary label and $Y_i^{(4)}\in\{1,2\}$. The final task is to identify if there are vertical ridges, which is also a binary label as $Y_i^{(5)}\in\{1, 2\}$.
%For convenience, we write $Y_i^{(k)}=\{1,\cdots, K^{(i)}\}$ with $(K^{(0)},\cdots, K^{(5)})=(45, 10, 5, 9, 2, 2)$, where $K^{(i)}$ is the total number of categories for the $i$-th task.}
See Figure \ref{fig:auxiliary-example} for the aforementioned auxiliary tasks of two specific products for visualization. The detailed categories of the five tasks can be found in Figure \ref{fig:auxiliary}.

\subsection{Model Training and Performance Measures}

We randomly split the whole dataset into a training set and a testing set, whose index sets are denoted by $\trainid$ and $\testid$, respectively.
To handle the photo information, we adopt the DenseNet architecture \citep{huang2017densely} for both the primary task and five auxiliary tasks.
First, we resize the input images to dimensions of 224 pixels for both height and width. Subsequently, we initialize the convolutional layers of the model using pre-trained weights provided by the open-source framework PyTorch \citep{paszke2017automatic}. After this pre-processing, the model should generate a 1920-dimensional feature vector after passing through the convolutional layers.
Then for the primary task, this feature vector undergoes mapping through a fully connected layer, which results in a total of 45 coefficient vectors. A softmax transformation is then applied to obtain the predicted probabilities.
Similarly, for each auxiliary task, the same feature vector is mapped through the corresponding fully connected layer with a softmax transformation.
To train the model, we employ the cross-entropy loss and utilize the Adam optimizer with a fixed learning rate of $10^{-4}$ so that the whole model parameters can be fine-tuned on our dataset. Similar operations are conducted for the auxiliary tasks.

After model training, we adopt the classification error rate (Err) to evaluate the prediction performance for the primary task on the testing set. Specifically, denote the sample size of the testing set as $n^*=\vert\testid\vert$.
Then the classification error rate is calculated as $\text{Err}=(n^*)^{-1}\sum_{i\in\testid} \mathbb I(Y_i^{(0)}\neq\wh Y_i^{(0)})$, where $\mathbb I(\cdot)$ is an indicator function. In this work, we adopt the Err metric primarily because, in real smart vending machine operations, the average revenue from each correct transaction settlement is quite limited, but the loss resulting from an incorrect transaction can be substantial. In case of a wrong transaction settlement, some customers may apply for a refund from the merchant, resulting in minimal loss for the merchant. However, other customers might directly file complaints with authoritative regulatory agencies for compensation, typically leading to a penalty of 10 times the original price of the goods involved in the transaction.
Note that the profit from a single correct transaction is no more than 10\% of the original price. This implies that the loss from one incorrect transaction settlement could require at least $10/0.1=100$ correct transactions to compensate. This makes the classification error rate the most important performance metric used in smart vending machine operations.

\subsection{Task and Rank Selection}

For a practical application, the auxiliary tasks and the rank of the coefficient matrix need to be selected. Intuitively, only those auxiliary tasks truly helpful for reducing the final classification error should be included. Those redundant auxiliary tasks should be excluded. Otherwise, the classification error rate can be inflated.
The rank of the coefficient matrix should also be carefully estimated. On one hand, the estimated rank $d$ should not be too small. Recall that the optimal weight is contained in $\mathbb W$. Once $d$ is under-estimated, the resulting weight might not be contained in $\mathbb W$. Hence the resulting weighted estimator would be seriously biased. On the other hand, the estimated rank cannot be too large either. Once $d$ is over-estimated, the resulting weight might not be optimal, and hence $\mathcal P(w^*)$ could be sub-optimal. Consequently, the resulting weighted estimator might suffer unnecessarily large variability.
Therefore, both the auxiliary tasks and the rank of the coefficient matrix should be carefully selected. To this end, we develop here a two-step cross-validation method for task and rank selection.

We start with the task selection. We temporarily assume the coefficient matrix $B$ is of full rank. 
This is a case with $d$ potentially over-estimated instead of under-estimated. Then the resulting weighted estimator should be at least unbiased as we mentioned previously.
Next, we start with the full task set $\mS_0=\{1,2,3,4,5\}$, which collects the indices of all auxiliary tasks. This represents the case with all auxiliary tasks being used for model training.
We next randomly split the whole dataset into $\trainid$ and $\testid$. Recall that $\trainid$ stands for the indices of the training dataset, which accounts for about $80\%$ of the whole dataset. Then the initial model can be computed based on $\mS_0$ and $\trainid$. The corresponding classification error rate is then evaluated on the $\testid$, which accounts for the rest 20\% of the whole dataset.  We replicate the experiment for a total of $50$ times by randomly partitioning the training and testing sets. This leads to a total of 50 Err values, which are then averaged and recorded as $\avgerr_0$.

For each $k\in\mS_0$, we consider a reduced task set as $\mS_{(0, k)}=\mS_0\backslash\{k\}$. We can then evaluate the classification error rate of $\mS_{(0, k)}$ by a similar method as for $\mS_0$. The classification error rates in 50 replications  are then averaged as $\avgerr_{(0, k)}$ for each $\mS_{(0, k)}$.
Next define $k_1^*=\arg\min_{k\in\mS_0} \avgerr_{(0, k)}$. This suggests that the auxiliary task indexed by $k_1^*$ should be excluded from $\mS_0$. This leads to an optimal choice of the auxiliary task set with size $(K-1)$ as $\mS_1=\mS_0\backslash\{k_1^*\}$. The corresponding $\avgerr$ value is then denoted by $\avgerr_1=\avgerr_{(0, k_1^*)}$.
We then repeat this process for every $\mS_k$ so that $\mS_{k+1}$ can be constructed for $0\leq k<5$, each of which is associated with $\avgerr_{k+1}$. This leads to a sequence of nested task sets $\mS_k$ with $0\leq k\leq 5$.
Then the best auxiliary task set is given by $\mS_{\text{opt}}=\mS_{k_{\text{opt}}}$, where $k_{\text{opt}}=\arg\min_k \avgerr_k$.
By applying this method to our smart-vending-machine dataset, we obtain the optimal auxiliary task set as $\mS_{\text{opt}}=\mS_1=\{1,2,3,4\}$ with $\avgerr_1=0.0468$.

\begin{figure}[h]
    \centering

    \subfloat{
        \includegraphics[width=.44\textwidth]{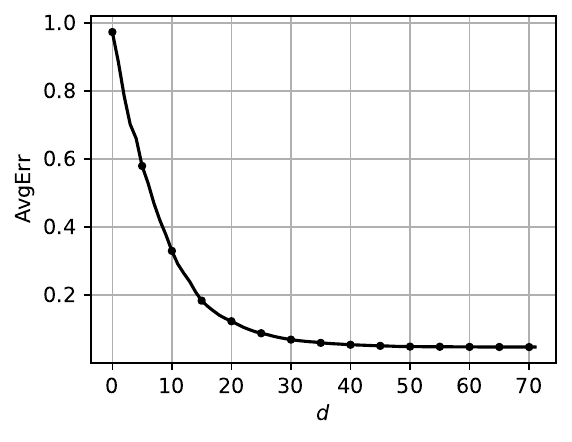}
    }
    \subfloat{
        \includegraphics[width=.47\textwidth]{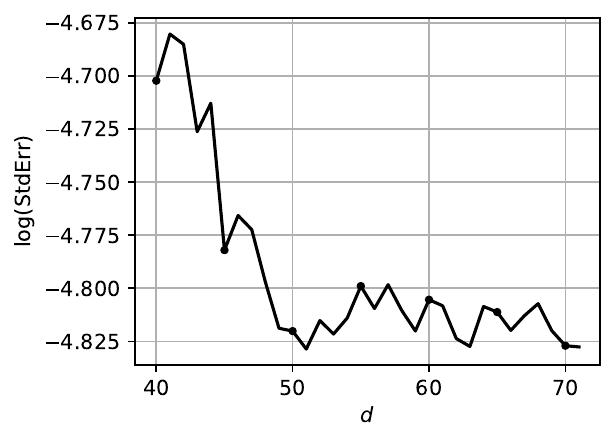}
    }

    \caption{The detailed results of AvgErr and StdErr (in logarithm) under different ranks.
    }
    \label{fig:rankselect}
\end{figure}

With the chosen $\mS_{\text{opt}}=\{1,2,3,4\}$, we next focus on the rank selection.  Specifically, we aim to investigate the impact of different values of $d$ on the classification error rate. Note that, with the auxiliary task set $\mS_{\text{opt}}=\{1,2,3,4\}$, the coefficient matrix $B$ is a $p\times K_0$ matrix, where $K_0=45+10+5+9+2=71$ is the total number of categories in the primary task and four auxiliary tasks. Therefore, we vary $d=\rank(B)$ in the range $[1, 71]$. For each given $d$, we compute the weighted estimator for the primary task on the training dataset $\trainid$. Then its classification error rate is evaluated on $\testid$. This process is also randomly replicated for a total of 50 times for each $d$. The resulting 50 Err values are then averaged as $\avgerr^{(d)}$ for each $d$. We then report $\avgerr^{(d)}$ in the left panel of Figure \ref{fig:rankselect}.

From Figure \ref{fig:rankselect}, we find a rather flat pattern for $\avgerr^{(d)}$ as $d\in [40, 71]$.
This suggests that the classification error rate is fairly insensitive for a broad range of $d\in[40, 71]$. Within this specification range of $d$, we also study the stability of the classification error rate by calculating the standard deviation of Errs, which is denoted as $\text{StdErr}^{(d)}$. We report this result in the right panel of Figure \ref{fig:rankselect}. We find that the minimum of $\text{StdErr}^{(d)}$ is achieved at $d=51$. Therefore, we should consider $d=51$ as the final specification about $\rank(B)$, as it yields almost the lowest classification error rate and the greatest stability.

\subsection{Model Comparison}

We compare the classification performance of our weighted estimator with that of other methods. With $\mS_{\text{opt}}=\{1,2,3,4\}$ and $d=51$, the final classification error rate of our method is given by $\avgerr=0.0481$. For comparison purpose, we consider a total number of five competitors. The first one is a benchmark model only using the primary task, which is denoted by STL (single-task learning). The second one is a model using all five auxiliary tasks, which can be considered as a special case of multi-task learning and thus denoted by MTL. Three auxiliary learning algorithms are also included for evaluation. They are, respectively, the AL-MTL method \citep{liebel2018auxiliary}, the AuxiLearn method \citep{navon2020auxiliary}, and the JTDS method \citep{chen2022auxiliary}. These three methods focus on obtaining the estimator by optimizing the weighted loss functions. For each method, the experiments are replicated for 50 times. The classification error rates of different methods are summarized in Table \ref{tab:comparison-dl}. The corresponding standard deviations ($\text{StdErr}$) for different methods are also reported. 
Additionally, we calculate the Relative Improvement (RI) in Errs between each competitor and our weighted estimator. Take the STL method as an example. The relative improvement is computed as $(1-0.0481/0.0535)\times 100\%=10.09\%$. Similarly, the RI value for the MTL method verse our method is computed as $(1-0.0481/0.0534)\times 100\%=9.93\%$.

\begin{table}[h]
\caption{The comparison results with different competitive methods. }
\centering
\begin{tabular}{cccc}
%S[table-format=+0.2,table-column-width=1.5cm]
\toprule
Methods   & AvgErr  & StdErr  & RI \\\midrule
Weighted Estimator (Ours) & 0.0481   & 0.0079  &$-$ \\
Single-Task Learning (STL)  & 0.0535 & 0.0211 & 10.09\% \\
Multi-Task Learning (MTL)  & 0.0534 & 0.0151 & 9.93\% \\
AL-MTL \citep{liebel2018auxiliary}    & 0.0808  & 0.0062  & 40.47\% \\
AuxiLearn \citep{navon2020auxiliary} & 0.0842   & 0.0251  & 42.87\% \\
JTDS \citep{chen2022auxiliary}     & 0.0766   & 0.0048  & 37.20\% \\  
\bottomrule             
\end{tabular}
\label{tab:comparison-dl}
\end{table}

From Table \ref{tab:comparison-dl}, we can draw the following conclusions. First, the classification error rates of three deep-learning based methods (i.e., AL-MTL, AuxiLearn, and JTDS) are much higher than those of the other methods. This finding suggests that, deep-learning methods are not applicable in our smart-vending-machine application.
Second, the classification error rate of our proposed weighted estimator outperforms that of STL and MTL. This result implies that, using carefully selected auxiliary tasks can benefit the prediction performance. Last, we find the standard deviation of our proposed method is much smaller than that of STL and MTL. This is consistent with our theoretical results in Corollary \ref{coro1}, since our method enjoys the optimal $\mathcal P(w^*)$ and thus has lower variability.

\section{CONCLUDING REMARKS}

Auxiliary learning is a popular method in machine learning research. The key idea of auxiliary learning is to introduce a set of auxiliary tasks to enhance the performance of the primary task. In this work, we adopt the idea of auxiliary learning to address the parameter estimation problem in high-dimensional settings. We start with the linear regression model. To improve the statistical efficiency of the parameter of primary research, we develop a weighted estimator, which is a linear combination of the OLS estimators of both the primary task and auxiliary tasks. The optimal weight is analytically derived and the statistical properties of the corresponding weighted estimator are studied. The weighted estimator is then extended to generalized linear regression models. Extensive numerical studies are presented to demonstrate the finite sample performance of the proposed method. Last, we apply the weighted estimator to address the product classification problem of smart vending machines. Results show that by using auxiliary learning, the products can be predicted more accurately than several state-of-the-art methods.

To conclude this article, we consider here some interesting topics for future studies.
First, we demonstrate that the proposed ORACLE and FEASIBLE estimators outperform the traditional OLS estimator by achieving lower mean squared errors. Our numerical studies suggest the usefulness of the weighted estimator in generalized linear models. However, the theoretical analysis under generalized linear models or more general loss functions needs further investigation.
Second, although we demonstrate that auxiliary learning can enhance estimation accuracy and efficiency, designing high-quality auxiliary tasks remains a significant challenge. This is an interesting topic for future exploration.

	%%%%%%%%%%%%%%%%%%%%%%%%%%%%%%%%%%%%%%%%%%%%%%%%%%%%%%%%%%%%%%%%%%%%%%%%%%%%%%%%%%%%%%%%%%%%%%%%%%%%%%%%%%%%%%%%%%%%%%%%%%%%
\section*{Acknowledgements}
	
This work is supported by National Natural Science Foundation of China (No.72371241, 72171229), Beijing Social Science Fund (24GLC033), and the MOE Project of Key Research Institute of Humanities and Social Sciences (22JJD910001).
	
	%%%%%%%%%%%%%%%%%%%%%%%%%%%%%%%%%%%%%%%%%%%%%%%%%%%%%%%%%%%%%%%%%%%%%%%%%%%%%%%%%%%%%%%%%%

\section*{Supplementary Materials}

Appendix A provides the detailed verification for the form of $\mathbb W$ and $w^*$. Appendix B provides the proof of Theorem \ref{thm-1}.
Appendix C provides some useful lemmas for Theorem \ref{thm-1}.

%\iffalse
\bibhang=1.7pc
\bibsep=2pt
\fontsize{9}{14pt plus.8pt minus .6pt}\selectfont
\renewcommand\bibname{\large \bf References}
%\begin{thebibliography}{11}
\expandafter\ifx\csname
natexlab\endcsname\relax\def\natexlab#1{#1}\fi
\expandafter\ifx\csname url\endcsname\relax
  \def\url#1{\texttt{#1}}\fi
\expandafter\ifx\csname urlprefix\endcsname\relax\def\urlprefix{URL}\fi
%\fi

\newpage
\bibliographystyle{chicago}      % Chicago style, author-year citations
\bibliography{ref}   % name your BibTeX data base

%-------------------------------------------
\vskip .65cm
\noindent
Guanghua School of Management, Peking University, Beijing, China.
\vskip 2pt
\noindent
E-mail: hanchao.yan@stu.pku.edu.cn
\vskip 2pt

\noindent
Center for Applied Statistics and School of Statistics, Renmin University of China, Beijing, China
\vskip 2pt
\noindent
E-mail: feifei.wang@ruc.edu.cn
\vskip 2pt

\noindent
Chinese Academy of International Trade and Economic Cooperation, Ministry of Commerce People's Republic of China, Beijing, China.
\vskip 2pt
\noindent
E-mail: xiachuanxin@caitec.org.cn
\vskip 2pt

\noindent
Guanghua School of Management, Peking University, Beijing, China.
\vskip 2pt
\noindent
E-mail: hansheng@pku.edu.cn
\vskip 2pt

\end{document}